\documentclass[reqno]{amsart}
\usepackage{amsmath}
\usepackage{amssymb}
\usepackage{stmaryrd}
\allowdisplaybreaks[2]
\title{
Linear Ramsey Numbers for Bounded-Degree Hypergrahps
}
\author{Yoshiyasu Ishigami}
\address{Department of Information and Communication Engineering, 
The University of Electro-Communications, Chofu, Tokyo 182-8585, Japan.
}
\email{yoshiyas@ice.uec.ac.jp}
\subjclass[2000]{05D55,05C65}
\keywords{Ramsey theory, Szemer\'edi's regularity lemma, hypergraph regularity lemma.}
\date{\today}
\def\picture #1 by #2 (#3){
                \vbox to #2{
                        \hrule width #1  height 0pt depth
0pt
                        \vfill
                        \special{picture #3}}}
\def\scaledpicture #1 by #2(#3 scaled #4){{
                \dimen0=#1 \dimen1=#2
                \divide\dimen0 by 1000 \multiply\dimen0 by
#4
                \divide\dimen1 by 1000 \multiply\dimen1 by
#4
                \picture\dimen0 by \dimen1 (#3 scaled #4)}}
\newtheorem{adf}{Definition}[section]
\newtheorem{tha}{Theorem}[section]

\newtheorem{thm}{Theorem}[section]
\newtheorem{cor}[thm]{Corollary}
\newtheorem{lemma}[thm]{Lemma}

\newtheorem{claim}[thm]{Claim}

\newtheorem{remark0}[thm]{Remark}

\newtheorem{algorithm0}[thm]{Algorithm}
\newtheorem{asett}[thm]{Setup}

\newenvironment{df}{
\begin{adf}\begin{sl}}
{\end{sl}\thqed\end{adf}}
\newenvironment{sett}{
\begin{asett}\begin{sl}}{\end{sl}\thqed\end{asett}}

\newcommand{\thqed}{\hfill\fbox{}\\ }

\def\brkt#1{\left({#1}\right)}

\def\ang#1{\langle{#1}\rangle}

\newcommand{\naturalset}{{\mathbb N}}

\newcommand{\Prob}{{\mathbb P}}
\newcommand{\Ex}{{\mathbb E}}

\newtheorem{procedure0}{Procedure}

\newcount\timehh\newcount\timemm\timehh=\time\divide\timehh by 60
\timemm=\time\count255=\timehh\multiply\count255 by-60
\advance\timemm by \count255
\ifnum\timehh=12\def\apm{pm}\else
\ifnum\timehh>12\def\apm{pm}\advance\timehh by-12\else
\def\apm{am}\fi\fi
\def\timestamp{\number\timehh\,:\,\ifnum\timemm<10 0\fi\number\timemm\,\apm}
\begin{document}
\maketitle
\begin{abstract}
We show that the the Ramsey number of 
every bounded-degree uniform hypergraph is linear 
with respect to the number of vertices.
This is a hypergraph extension of the famous theorem for ordinary graphs which 
Chv\'atal et al. \cite{CRST} showed in 1983.
Our result may demonstrate the potential of a new 
hypergraph regularity lemma by \cite{I06}.
\end{abstract}
\section{Introduction}
A {\bf $k$-uniform hypergraph} is a family of $k$-element subsets (called 
\lq edges\rq) of the underlying set, whose members are called \lq vertices.\rq\ 
It is {\bf complete} if and only if it contains all the $k$-element subsets.
For a $k$-uniform hypergraph $H$ and a positive integer $b$,
 the {\bf Ramsey number of } $H$, denoted by $R_b(H)$, is 
the least integer $R$ such that for any $b$-coloring of 
the edges of the $k$-uniform complete hypergraph on $R$ vertices, 
there exists a monochromatic copy of $H$. 
The study of this number is a main theme of 
Ramsey Theory, which 
has been considered to be a central field of 
 combinatorics or discrete mathematics.
Ramsey theory started by the following theorem. 
\begin{tha}[Ramsey (1930) \cite{Ramsey}]\label{xx17}
Let $b,k$ and $N$ be positive integers.
For any $k$-uniform hypergraph $H$ on $N$ vertices,
its Ramsey number $R_b(H)$ exists.
\end{tha}
As one of the earliest deep applications of the regularity lemma by 
Szemer\'edi, the following fundamental theorem in Ramsey theory was obtained.
It was a conjecture of Burr and Erd\"{o}s \cite{BE}.
\begin{tha}[Chv\'atal-R\"{o}dl-Szemer\'edi-Trotter (1983)\cite{CRST}]
\label{xx18a1}
Let $b\ge 1$ be a constant integer and $N\ge 1$ be a (large) integer.
For any ordinary graph (i.e. a $2$-uniform hypergraph) on 
$N$ vertices with maximum degree $O(1)$, we have 
 $R_b(H)={O_b(N)}.$
\end{tha}
For a hypergraph, we say that a vertex is a {\bf neighbor} of another different vertex 
if-and-only-if there exists an edge containing the two vertices.
The {\bf degree} of a vertex is the number of neighbors of the vertex.
The {\bf maximum degree} of a hypergraph is defined to be 
the largest degree over all vertices.\par
Very recently, two groups obtained the following independently by different methods, 
though both depend on the hypergraph regularity platform 
of Frankl-R\"{o}dl (2002) \cite{FR02}.
\begin{tha}[Cooley et al.\cite{CFKO} and Nagle et al.\cite{NORS}]
\label{xx18a2}
Let $N$ be a (large) integer.
For any $3$-uniform hypergraph on $N$ vertices with maximum degree $O(1)$,
 we have  $R_2(H)=O(N)$.
\end{tha}
Kostochka-R\"{o}dl (2006) \cite{KR} showed that 
$R_2(H)\le N^{1+o(1)}$ for any $O(1)$-uniform hypergraph on $N$ vertices with maximum degree 
$O(1).$
In this paper, we will prove the following theorem.
\begin{thm}[Main Theorem]\label{xx16x}
Let $b\ge 1$ be a constant integer and $N\ge 1$ be a (large) integer.
For any $O(1)$-uniform hypergraph on $N$ vertices with maximum degree $O(1)$,
 we have $R_b(H)=O_{b}(N).$
\end{thm}
\bigskip
\par
I uploaded the first draft \cite{I06lr} of this result to 
the preprint server, arxiv.org 
(http://arxiv.org/), on 20 Dec. 2006.
After writing almost all parts of it, I learned the existence of 
a preprint by 
Cooley et al. \cite{CFKO2} uploaded to the preprint server 
on 13 Dec. 2006. They 
obtained the two-color case of the 
main theorem independently from us.
However, our method is different from theirs.
Their method relies on a regularity lemma with a counting lemma by 
R\"{o}dl-Schacht \cite{RS}, which need long proofs.
(\cite{RS} is not self-contained. It uses results from \cite[Th.6.5,Cor.6.11]{KNR} 
and omits technical proofs (\cite[Prop.28,29,30,32,33]{RS}) which are 
straightforward or similar to proofs in \cite{FR02,NRS,RSk04}.)
On the other hand, the version of the regularity lemma from \cite{I06} which we will use
 has a short proof.
While our proof is simple, the main lemma(Lemma \ref{xx11} or Corollary \ref{071020}, 
counting lemma for blowups) 
is stronger than their corresponding main lemma(they called the embedding lemma), 
since our regularity setting is weaker in a sense.
\par
The main purpose of this paper is not only to prove the fundamental theorem 
in Ramsey theory but also to show the potential of the framework of \cite{I06}.
Although another proof of Theorem \ref{xx18a1} without the graph 
regularity lemma \cite{Sz} was found later in \cite{GRR}, 
the techniques developed in \cite{CRST} have been used 
for many applications. It may be why Theorem \ref{xx18a1} is considered 
as a milestone in the survey \cite{KSSS02}. 
\cite[\S5.1]{KSSS02} says that \cite{CRST} was probably the first deep application 
of the regularity lemma. (On the other hand, 
Chvatal-Szemer\'edi \cite{ChSz} was published earlier and also deep, and 
some techniques of \cite{CRST} appeared already in \cite{ChSz}.
The main theorem in \cite{ChSz} is extended in \cite{I}.) 
I believe that the technique of this paper will be used for other applications.
Such an example can be seen already in \cite{I06m}.
\par
The regularity lemma by \cite{I06} gives a new proof of the 
Szemeredi theorem on progressions 
which is shorter than previous proofs.
Due to the simplicity of the proof, it is not hard to modify 
the proof of the regularity lemma for deeper applications if necessary.
Although we need only the surface of the theorem for the purpose of this paper,
we already have an application which needs 
 a slight modification of our regularity lemma. 
See \cite{I06} for discussion on differences from earlier hypergraph regularity 
lemmas \cite{RSk04,NRS,G,Tao06,RS}.
\par
Cooley et al. \cite{CFKO, CFKO2} and Nagle et al. \cite{NORS} treated 
only the case of $2$-coloring. 
Although their methods may be essentially extendable to the multicolor case, 
it should need more technical work and pages in their setting.
On the other hand, from the beginning plan of our regularity lemma, 
we have considered the multicolor case because it is natural for both of regularity lemma and 
its applications.
\section{Statements of Regularity Lemma and Main Lemma}
In this paper, we denote by $\Prob$ and $\Ex$ the probability and expectation, 
respectively. We denote the conditional probability and exepctation by
$\Prob[\cdots|\cdots]$ and $\Ex[\cdots|\cdots].$ 
\begin{sett}\label{r051231}
Throughout this paper, we fix a positive integer $r$ and 
an \lq index\rq\ set $\mathfrak{r}$ with $|\mathfrak{r}|=r.$ 
Also we fix a probability space 
$({\bf \Omega}_i,{\mathcal B}_i,\Prob)$ 
for each $i\in \mathfrak{r}$.
Assume that ${\bf \Omega}_i$ is finite (but its cardinality may not be 
constant) 
and ${\mathcal B}_i=2^{{\bf \Omega}_i}$ 
for the sake of simplicity.
Write ${\bf \Omega}:=({\bf \Omega}_i)_{i\in \mathfrak{r}}$.
\end{sett}
In order to avoid using technical words like mesurability or 
Fubini's theorem frequently to readers who are interested only in applications to 
discrete mathematics, 
we assume ${\bf \Omega}_i$ as a (non-empty) finite set.
However our argument should be extendable to a more general probability space.
For applications, ${\bf \Omega}_i$ would contain a huge number of vertices, though 
we do not use the assumption in our proof.
\par
For an integer $a$, we write $[a]:=\{1,2,\cdots,a\},$ and 
${\mathfrak{r}\choose [a]}:=\dot{\bigcup}_{i\in [a]}{\mathfrak{r}\choose i}
=\dot{\bigcup}_{i\in [a]}\{I\subset \mathfrak{r}| |I|=i\}.$
When $r$ sets $X_i, i\in {\mathfrak r},$ with indices from ${\mathfrak r}$ are 
called {\bf vertex sets}, 
we write $X_J:=\{e\subset \dot{\bigcup}_{i\in J}X_i| |e\cap X_j|=1 \forall j\in J\}$ 
whenever $J\subset {\mathfrak r}$.
\begin{df}[(Colored hyper)graphs]
Suppose Setup \ref{r051231}. 
A {\bf $k$-bound $(b_i)_{i\in [k]}$-colored ($\mathfrak{r}$-partite hyper)graph} $H$ 
is a triple $((X_i)_{i\in\mathfrak{r}},({C}_I)_{I\in {\mathfrak{r}\choose [k]}},
(\gamma_I)_{I\in {\mathfrak{r}\choose [k]}}
)$ where (1) each $X_i$ is a set called a \lq vertex set,\rq\ (2)
${C}_I$ is a set with at most $b_{|I|}$ elements, and 
(3) $\gamma_I$ is a function from $X_I$ to ${C}_I.$
We write $V(H)=\dot{\bigcup}_{i\in \mathfrak{r}}
V_i(H)=\dot{\bigcup}_{i\in \mathfrak{r}}X_i
$ and 
${\rm C}_I(H)={C}_I.$ 
Each element of $V(H)$ is called a {\bf vertex}.
Each element $e\in V_I(H)=X_I, I\in {\mathfrak{r}\choose [k]},$ is called 
 an {\bf (index-$I$ size-$|I|$) edge}.
Each member in ${\rm C}_I(H)$ is a {\bf (face-)color (of index $I$)}.
Write $H(e)=\gamma_I(e)$ for each $I.$ 
\par
Let $I\in {\mathfrak{r}\choose [k]}$ 
and $e\in V_I(H).$ 
For another index $\emptyset\not=
J\subset I$, we denote by $e|_J$ the index-$J$ edge 
$e\setminus \brkt{\bigcup_{j\in I\setminus J}X_j}\in V_J(H)$.
We define the {\bf frame-color}
and {\bf total-color} of $e$ by ${H}(\partial
e):=({H}(e|_J)|\,\emptyset\not=J\subsetneq I)$ and by
${H}(\ang{e})=H\ang{e}:=({H}(e|_J)|\,\emptyset\not=J\subsetneq I).$
Write ${\rm TC}_I(H):=\{H\ang{e}|\,{e}\in X_I\},$ ${\rm
TC}_s(H):=\bigcup_{I\in {\mathfrak{r}\choose s}}{\rm TC}_I(H),$
and ${\rm TC}(H):=\bigcup_{s\in  [k]}{\rm TC}_s(H).$
\par
A {\bf ($k$-bound) (simplicial-)complex} is a $k$-bound
(colored ${\mathfrak r}$-partite hyper)graph 
such that for each $I\in {\mathfrak{r}\choose [k]}$ 
there exists at most one index-$I$ color called \lq invisible\rq\ 
and that if (the color of) an edge $e$ is invisible then 
any edge $e^*\supset e$ is invisible. An edge or its color 
is {\bf visible} if it is not invisible. 
\par
For a $k$-bound graph ${\bf G}$ on ${\bf \Omega}$ and $s\le k$, 
let ${\mathcal S}_{r,s,h,{\bf G}}={\mathcal S}_{s,h,{\bf G}}$ be the set of $s$-bound 
simplicial-complexes $S$ such that 
(1) each of the $r$ vertex sets contains exactly $h$ vertices and 
that (2)
for any $I\in {\mathfrak{r}\choose [s]}$ 
 there is an injection from the index-$I$ visible colors of $S$ to the 
index-$I$ colors of ${\bf G}$.
(When a visible color $\mathfrak{c}$ of $S$ corresponds to another color $\mathfrak{c}'$ of 
${\bf G}$, we simply write $\mathfrak{c}=\mathfrak{c}'$ 
without presenting the injection explicitly.)
For $S\in {\mathcal S}_{s,h,{\bf G}}$, we denote by ${\mathbb V}_I(S)$ the set of 
index-$I$ visible edges. Write ${\mathbb V}_i(S):=\bigcup_{I\in {\mathfrak{r}\choose 
i}}{\mathbb V}_I(S)$ and ${\mathbb V}(S):=\bigcup_i {\mathbb V}_i(S).$
\par
For a complex $S$ and $U\subset V(S)$, we denote by 
$S\setminus U$ the complex obtained from $S$ by deleting the vertices in $U$ and 
the edges containing a vertex in $U.$ When $U$ consists of a single vertex $u$, 
we write $S\setminus \{u\}=S\setminus u$.
Also write $S\setminus V(N)=S\setminus N$ for another complex $N.$
Sometimes we write $S|_U=S\setminus(V(S)\setminus U)$ and call it 
the {\bf complex of $S$ induced by $U$}.
\end{df}
\begin{df}[Partitionwise maps]
A {\bf partitionwise map} $\varphi$ is 
a map 
 from $r$ vertex sets $W_i,i\in \mathfrak{r},$ with $|W_i|<\infty$ to {\em the} 
$r$ vertex sets (probability spaces)$
U_i,i\in\mathfrak{r}$,
such that 
each $w\in W_i$ is mappped into $U_i.$
We denote by $\Phi((W_i)_{i\in\mathfrak{r}}, (U_i)_{i\in\mathfrak{r}})$ 
or $\Phi(\bigcup_{i\in\mathfrak{r}}W_i,\bigcup_{i\in\mathfrak{r}}U_i)$ 
the set of partitionwise maps from $(W_i)_i$ to $(U_i)_i$.
If $U_i={\bf \Omega}_i$ or $U_i$ is obvious then we omit them.
A partitionwise map is {\bf random} if and only if 
each $w\in W_i$ is mutually-independently mapped at random 
according to the probability space 
${\bf \Omega}_i$.
\par
For two partitionwise maps $\phi\in \Phi((W_i)_i)$ and $\phi'\in \Phi((W'_i)_i),$ 
we denote by $\phi\dot{\cup} \phi'$ 
the partitionwise map $\phi^*\in \Phi((W_i\dot{\cup} W'_i)_i)$ such that 
$\phi^*(w)=\phi(w)$ and $\phi^*(w')=\phi'(w')$ for all $w\in W_i,w'\in W'_i, i\in
\mathfrak{r}$, where if $W_i\cap W'_i\not=\emptyset$ for some $i$
 then we consider a copy of $W'_i$ so that the two domains are disjoint.
\par
Sometimes for a graph (a complex, usually) $S$, 
we write $\Phi(V(S))=\Phi(S)$ when it is not confusing.
\par
For two $r$-partite graphs $S,G$ and for a partitionwise map $\phi\in \Phi(W,V(G))$ 
with some $W\supset V(S)$, 
we say that $\phi$ {\bf embeds} $S$ in $G$, or write 
\begin{align}
S\stackrel{\phi}{\hookrightarrow}G
\end{align}
if and only if $S(e)=G(\phi(e))$ for all $e\in \mathbb{V}(S).$ 
\end{df}
Suppose that $\phi$ is random and that any two 
events $S(e)=G(\phi(e))$ and $S(e')=G(\phi(e'))$ are mutually 
independent 
unless $e=e'$.
(This happens if all edges of $G$ are colored uniformly at random.)
Then 
we observe that 
\begin{eqnarray*}
\Prob_{\phi\in\Phi(S)}[S\stackrel{\phi}{\hookrightarrow}G]
&=&\Prob_{\phi\in\Phi(S)}[
G(\phi(e))=S(e)
 \forall e\in\mathbb{V}(S)
]\\
&=&\prod_{I\in {\mathfrak{r}\choose [k]}}\prod_{e\in\mathbb{V}_I(S)}
\Prob_{\phi\in\Phi(S)}[G(\phi(e))=S(e)| G(\phi(e^*))
=S(e^*) \forall e^*\subsetneq e
]\\
&=&\prod_I\prod_{e\in\mathbb{V}_I(S)}
\Prob_{{\bf e}\in V_I(G)}[G({\bf e})=S(e)| G({\bf e}|_J)
=S(e|_J) \forall J\subsetneq I
]
\end{eqnarray*}
where ${\bf e}|_J$ and $e|_J$ are the edges restricted in index $J$.
With this observation, we define the regularity of hypergraphs.
\begin{df}[Regularity]
Let ${\bf G}$ be a $k$-bound graph on $
{\bf \Omega}$.
For 
$\vec{\mathfrak{c}}=(\mathfrak{c}_J)_{J\subset I}\in {\rm TC}_I({\bf G}), 
I\in {\mathfrak{r}\choose [k]}$, we define {\bf relative density}
\begin{eqnarray*}
{\bf d}_{\bf G}(\vec{\mathfrak{c}}):=
\Prob_{{\bf e}\in {\bf \Omega}_I
}[
{\bf G}({\bf e})=\mathfrak{c}_I
|
{\bf G}(\partial{\bf e})=
(\mathfrak{c}_J)_{J\subsetneq I}
].
\end{eqnarray*}
\par
For a positive integer $h$ and a function 
$\varepsilon:[k]\times \naturalset\to (0,1]$, we say that 
 ${\bf G}$ is 
 {\bf $(\varepsilon,h)$-regular} 
if and only if 
there exists a function
${ \delta}: {\rm TC}({\bf G})\to [0,\infty)$ 
such that 
\begin{eqnarray}
\hspace{-5mm}
{\rm (i)}&
 \Prob_{\phi\in\Phi(S)}
[
S\stackrel{\phi}{\hookrightarrow}{\bf G}
]=
\displaystyle\prod_{e\in {\mathbb V}(S)}
\brkt{
{\bf d}_{\bf G}(S\ang{e})
\dot{\pm}
\delta(S\ang{e})
}&
\forall S\in {\mathcal S}_{k,h,{\bf G}},
\nonumber
\\
\hspace{-5mm}
{\rm (ii)}& \Ex_{{\bf e}\in {\bf \Omega}_I}[\delta({\bf G}\ang{\bf e})]\le 
\varepsilon
\brkt{|I|,
\max_{J\in {\mathfrak{r}\choose [k]}:|J|\ge |I|
}
|{\rm C}_J({\bf G})|
}
& \forall I\in {\mathfrak{r}\choose [k]},
\nonumber
\end{eqnarray}
where $a\dot{\pm}b$ means (the interval of) numbers $c$ with 
$\max\{0,a-b\}\le c\le \min\{1,a+b\}$.
\par
A {\bf subdivision }  of a $k$-bound graph ${\bf G}$ on ${\bf \Omega}$ 
is a $k$-bound graph ${\bf G}^*$ on the same ${\bf \Omega}$ such that 
\\
(i) for any size-$k$ edge ${\bf e}\in {\bf \Omega}_I$ with 
$I\in {\mathfrak{r}\choose k},$ 
it holds that ${\bf G}^*({\bf e})={\bf G}({\bf e}),$ and \\
(ii) for any two edges ${\bf e}, {\bf e'}\in {\bf \Omega}_I$ with 
$I\in {\mathfrak{r}\choose [k-1]}$, 
if ${\bf G}^*({\bf e})={\bf G}^*({\bf e'})$ then ${\bf G}({\bf e})={\bf G}({\bf e'}).$
\end{df}
\begin{tha}[Hypergraph Regularity Lemma in \cite{I06}] 
\label{r060721}
Let $r\ge k,h,\vec{b}=(b_i)_{i\in [k]}$ be positive integers, and 
$\varepsilon:[k]\times \naturalset \to (0,1]$ a function.
Then there exist integers 
$\widetilde{b}_1\ge \cdots\ge \widetilde{b}_{k-1}
$ such that if ${\bf G}$  is 
a $\vec{b}$-colored ($k$-bound $r$-partite hyper)graph on ${\bf \Omega}$ then 
there exists an $(\varepsilon,h)$-regular $(\widetilde{b}_1,\cdots,
\widetilde{b}_{k-1},b_k)$-colored 
subdivision ${\bf G}^*$ of ${\bf G}.$ 
\end{tha}
\hspace{5mm} 
Two earliest versions of the hypergraph regularity lemmas 
were obtained by R\"{o}dl and his collaborators \cite{RSk04,RS} and by
Gowers \cite{G} independently, and another one was obtained by Tao \cite{Tao06}.
R\"{o}dl-Schacht \cite{RS} obtained a variant of their earlier one so that 
it would be more appropriate for applications.
(We discuss the differences between these regularity lemmas in \cite{I06}.)
(For earlier results about (weaker) 
hypergraph regularity lemmas, see 
\cite{Chung,Chung91,CG,CG2,CG3,CT,HT,HT2}.
)
\par
Theorem \ref{r060721} lacks an important part of the main theorem in \cite{I06}, 
the simple way to construct the subdivision. Although it is very important, 
we will not need it for our purpose of this paper.
\begin{df}[Blowup]
For a positive integer $\Delta$,
 a {\bf $\Delta$-blowup} of a complex $S$ 
is an 
($\mathfrak{r}$-partite) $k$-bound complex $B$ 
on a finite set of vertices 
with maximum degree $\Delta(B)\le \Delta$ 
such that 
 $B$ is embeddable in $S$ (i.e. 
 there exists a map $\phi$ which embeds $B$ in $S$)
where the {\bf maximum degree} of $B$ is defined by 
$$
\Delta(B):=\max_{v\in V(B)}
|\{w\in V(B)\setminus \{v\}
| \{v,w\}\in \mathbb{V}_2(B)
\}|.$$
Note that $\max_{v\in V(B)}
|\{e\in \mathbb{V}_k(B)|v\in e
\}|\le {\Delta(B)\choose k-1}.$
\end{df}
Our main theorem will be obtained as a collorary of 
Theorem \ref{r060721} and the following.
\begin{lemma}[Main Lemma - Counting Lemma for Blowups]\label{xx11}
For any positive integers 
$k$ and $\Delta,$ 
 there exist $k$ functions 
$\eta_i=\eta_i(\rho_i)>0, i\in [k],$
(independent from ${\bf \Omega}$) such that the following holds 
for any reals $
0<\rho_1\le \rho_2\le\cdots\le \rho_k< 1$.
\par
Let $r\ge k$ and $h$ be positive integers.
Let ${\bf G}$ be an ($\mathfrak{r}$-partite) 
$k$-bound ($(1/\rho_i)_i$-colored) hypergraph on (any probability space)
${\bf \Omega}=({\bf \Omega}_i)_{i\in\mathfrak{r}}.$
Let $S\in {\mathcal S}_{r,k,h,{\bf G}}.$
Suppose,  for any $(2\Delta)$-blowup $S'$ of $S$ with 
$|V(S')|\le 2\Delta^{2k}$, the property that
\begin{eqnarray}
\Prob_{\phi\in \Phi(V(S'))}[
S'\stackrel{\phi}{\hookrightarrow} {\bf G}
]
=
\prod_{e\in \mathbb{V}(S')}
(1\dot{\pm}\eta_{|e|}(\rho_{|e|}))
{\bf d}_{\bf G}(S'\ang{e})
\label{xx11a}
\end{eqnarray}
and further suppose that
\begin{eqnarray}
{\bf d}_{\bf G}(S\ang{e})>
\rho_{|e|}
\quad \forall e\in \mathbb{V}(S). \label{xx15d}
\end{eqnarray}
Let $B$ be a $\Delta$-blowup of $S$.
Then for any vertex $u\in V(B),$
\begin{eqnarray}
&&
\Ex_{\varphi\in\Phi(B\setminus u)}\left[\left.
\Prob_{\phi\in \Phi(\{u\})}[
B\stackrel{\varphi\dot{\cup}\phi}{\hookrightarrow} {\bf G}
]
\right|
B\setminus u\stackrel{\varphi}{\hookrightarrow} {\bf G}
\right]
\ge (1-\eta_k^{1/4})
\prod_{e\in \mathbb{V}(B):u\in e}
{\bf d}_{\bf G}(B\ang{e}).\label{xx16}
\end{eqnarray}
\end{lemma}
Of course, in the above, the exact value $1/4$ of (\ref{xx16}) is not important
 here. 
Note that each $\eta_i$ is independent from $\rho_j,j<i,$ and $|V(B)|<\infty$.
\begin{cor}\label{071020}
In Lemma \ref{xx11}, if each ${\bf \Omega}_i$ is a finite set and 
if $|V_i(B)|<\eta_1(\rho_1)|{\bf \Omega}_i|$ for each $i\in \mathfrak{r}$ then 
the left hand side of (\ref{xx16}) can be replaced by 
$$\Ex_{\varphi\in\Phi(B\setminus u)}\left[\left.
\Prob_{\phi\in \Phi(\{u\})}[
B\stackrel{\varphi\dot{\cup}\phi}{\hookrightarrow} {\bf G}
\mbox{ and } \phi(u)\not\in \varphi(V(B\setminus u))]
\right|
B\setminus u\stackrel{\varphi}{\hookrightarrow} {\bf G}
\right].$$
In particular, 
$$\Prob_{\varphi\in\Phi(B)}\left[
B\stackrel{\varphi}{\hookrightarrow} {\bf G}
\mbox{ and $\varphi$ is an injection} 
\right]
\ge (1-\eta_k^{1/4})^{|V(B)|}
\prod_{e\in\mathbb{V}(B)}{\bf d}_{\bf G}(B\ang{e})>0
.$$
\end{cor}
\section{Proof of Main Lemma}
Our proof concept is to repeat $k-1$ times of an argument which 
 Cooley et al. \cite{CFKO} repeated twice for the $3$-uniform case.
Cooley et al. \cite{CFKO2} avoided the iteration and employed the \lq half\rq\ dense 
version of the regularity lemma with the counting lemma by R\"{o}dl-Schacht \cite{RS}.
However the iteration will work smoothly in the platform of the regularity lemma by 
\cite{I06}.
\begin{df}[Abbreviation]
$\bullet$ We write \lq iff\rq\ for \lq if and only if\rq.\\
$\bullet$
For a complex $B$ and its edge $e$, we write 
${\bf d}^{(+\delta)}_{\bf G}(B\ang{e}):=\min\{1,
(1+\eta_{|e|})
{\bf d}_{\bf G}(B\ang{e})\},$
${\bf d}^{(-\delta)}_{\bf G}(\mathfrak{c}):=\max\{0,(1-\eta_{|e|})
{\bf d}_{\bf G}(B\ang{e})\},$
and 
${\bf d}^{(\delta)}_{\bf G}(B\ang{e}):=(1\dot{\pm}\eta_{|e|})
{\bf d}_{\bf G}(B\ang{e}).$
\\
$\bullet$ For a $k$-bound complex $B$ and an integer $i\le k$, 
we denote by $B^{\ang{i}}$ the complex obtained from $B$ by 
invisualizing all edges of size at least $i+1.$
That is, 
$\mathbb{V}_j(B^{\ang{i}})=\mathbb{V}_j(B)$ for all $j\le i$ and 
$\mathbb{V}_j(B^{\ang{i}})=\emptyset$ for all $j>i$.
\\
$\bullet$ For  a $k$-bound complex $B$ and an integer $i\le k$, 
write $\mathbb{V}_{(i)}(B):=\bigcup_{j\le i}\mathbb{V}_j(B).$
\\
$\bullet$
A complex $S'$ is a {\bf subcomplex} of another complex $S$ iff 
there exists an injection which embeds $S'$ in $S$.
\end{df}
We will prove Lemma \ref{xx11} by induction on $k$ and on $|V(B)|$.
If $k=1$ then it is trivial.
We assume that $k\ge 2$ and the assertion holds for $k-1$ or less, 
since $B$ has no edge of size $k$ in those cases.
When $|V(B)|<k$ then it is clear from the induction hypothesis.
Assume that $|V(B)|\ge k.$
\begin{df}
Let $B',B''$ be subcomplexes of $B$ such that $B'=B''|_{V(B')}.$
Let $\varphi\in\Phi(B')$ (or $\varphi\in\Phi(W)$ for some 
$W\supset V(B')$ with $(V(B'')\setminus V(B'))\cap W=\emptyset$).
If  $B'\stackrel{\varphi}{\hookrightarrow}{\bf G}$ then 
we define the {\bf extension error of $\varphi$ from $B'$ to $B''$} by
\begin{eqnarray}
\beta(\varphi,B',B''):=\left|
\Prob_{\phi\in\Phi(B''\setminus V(B'))}
[B''\stackrel{\varphi\dot{\cup}\phi}{\hookrightarrow}{\bf G}]
\brkt{
\prod_{e\in \mathbb{V}(B'')\setminus \mathbb{V}(B')}{\bf d}_{\bf G}(B\ang{e})}^{-1}
-1\right|^2.\label{xx15e}
\end{eqnarray}
When $B'$ is empty(or when all visible edges contain no vertex in $B'$), 
we can naturally define $\beta(B'')
=\left|
\Prob_{\phi\in\Phi(B'')}
[B''\stackrel{\phi}{\hookrightarrow}{\bf G}]
\brkt{
\prod_{e\in \mathbb{V}(B'')}{\bf d}_{\bf G}(B\ang{e})}^{-1}
-1\right|^2.
$
\end{df}
\begin{claim}[Extension error is usually small]\label{xx12b}
Let $\ell\le k$ and $B', B''$ be $\ell$-bound subcomplexes of 
$B$. Suppose that $B'=B''|_{V(B')}$ and 
that $|V(B')|\le |V(B'')|\le \Delta^{2k}$. (Without loss of generality, 
$\Delta\ge 2.$)
Then we see that 
\begin{eqnarray*}
\Ex_{\varphi\in\Phi(B')}\left[\left.
\beta(\varphi,B',B'')
\right|
B'
\stackrel{\varphi}{\hookrightarrow} {\bf G}
\right]
\le \beta_\ell
\end{eqnarray*}
where $\beta_\ell:=
{8.1\cdot 2^{\Delta^{2k}} \eta_{\ell}}$.
\end{claim}
\begin{proof} 
We consider the following complex $B^*$.
We let $V(B^*)=V(B')\dot{\cup}W\dot{\cup}W'$ 
where $W,W'$ are two disjoint copies of $V(B'')\setminus V(B').$
Every edge $e$ with $|e\cap W|\cdot |e\cap W'|>0$
is invisible in $B^*.$
Any other edges have the same colors as the corresponding edges of $B''.$
(That is, $B^*$ is obtained from $B''$ by \lq splitting\rq\ 
 $B''\setminus V(B').$)
Since $|V(B^*)|\le 2\Delta^{2k},$
the assumption (\ref{xx11a}) yields the property that 
\begin{eqnarray*}
&&
\Ex_{\varphi\in \Phi(B')}
[
\brkt{
\Prob_{\phi\in \Phi(B''\setminus V(B'))}
[B''\stackrel{\varphi\dot{\cup}\phi}{\hookrightarrow}
{\bf G}]
-\prod_{e\in\mathbb{V}(B'')\setminus\mathbb{V}(B')}{\bf d}_{\bf G}(B\ang{e})
}^2
|B'
\stackrel{\varphi}{\hookrightarrow}{\bf G}
]
\\
&=&
\Ex_{\varphi\in \Phi(B')}
[
\brkt{
\Prob_{\phi\in \Phi(B''\setminus V(B'))}
[B''\stackrel{\varphi\dot{\cup}\phi}{\hookrightarrow}
{\bf G}]
}^2
|B'\stackrel{\varphi}{\hookrightarrow}{\bf G}
]
+\brkt{
\prod_{e\in\mathbb{V}(B'')\setminus\mathbb{V}(B')}{\bf d}_{\bf G}(B\ang{e})
}^2
\\
&&
-2
\prod_{e\in\mathbb{V}(B'')\setminus\mathbb{V}(B')}{\bf d}_{\bf G}(B\ang{e})\cdot
\Prob_{\phi\in \Phi(B'')}
[
B''\stackrel{\phi}{\hookrightarrow}
{\bf G}
|B'\stackrel{\phi}{\hookrightarrow}{\bf G}
]
\\
&\stackrel{(\ref{xx11a})}{=}&
\Prob_{\phi\in \Phi(B^*)}
[B^*\stackrel{\phi}{\hookrightarrow}
{\bf G}
|B'\stackrel{\phi}{\hookrightarrow}{\bf G}
]
+\brkt{
\prod_{e\in\mathbb{V}(B'')\setminus\mathbb{V}(B')}{\bf d}_{\bf G}(B\ang{e})
}^2
\\
&&
-2\brkt{
\prod_{e\in\mathbb{V}(B'')\setminus\mathbb{V}(B')}{\bf d}_{\bf G}(B\ang{e})
}\prod_{e\in\mathbb{V}(B'')}{\bf d}_{\bf G}^{(\delta)}(B\ang{e})
/
\prod_{e\in\mathbb{V}(B')}
{\bf d}_{\bf G}^{(\delta)}(B\ang{e})
\\
&\stackrel{(\ref{xx11a})}{=}&
\prod_{e\in \mathbb{V}(B^*)}
{\bf d}^{(\delta)}_{\bf G}(B^*\ang{e})
/\prod_{e\in \mathbb{V}(B')}
{\bf d}^{(\delta)}_{\bf G}(B\ang{e})
+\brkt{
\prod_{e\in\mathbb{V}(B'')\setminus\mathbb{V}(B')}{\bf d}_{\bf G}^{(\delta)}(B\ang{e})
}^2
\\
&&
-2\brkt{
\prod_{e\in\mathbb{V}(B'')\setminus\mathbb{V}(B')}{\bf d}_{\bf G}^{(\delta)}(B\ang{e})
}^2\prod_{e\in\mathbb{V}(B')}
{
{\bf d}_{\bf G}^{(\delta)}(B\ang{e})
\over
{\bf d}_{\bf G}^{(\delta)}(B\ang{e})
}
\\
&\le&
\brkt{
\prod_{e\in\mathbb{V}(B'')\setminus\mathbb{V}(B')}{\bf d}_{\bf G}^{(+\delta)}(B\ang{e})
}^2
\brkt{\brkt{
\prod_{e\in\mathbb{V}(B')}
{
{\bf d}_{\bf G}^{(+\delta)}(B\ang{e})
\over
{\bf d}_{\bf G}^{(-\delta)}(B\ang{e})
}}+1
}
\\
&&
-2\brkt{
\prod_{e\in\mathbb{V}(B'')\setminus\mathbb{V}(B')}{\bf d}_{\bf G}^{(-\delta)}(B\ang{e})
}^2\prod_{e\in\mathbb{V}(B')}
{
{\bf d}_{\bf G}^{(-\delta)}(B\ang{e})
\over
{\bf d}_{\bf G}^{(+\delta)}(B\ang{e})
}
\\
&\le&
\brkt{
\prod_{e\in\mathbb{V}(B'')\setminus\mathbb{V}(B')}
{\bf d}_{\bf G}(B\ang{e})
}^2
\\
&&\cdot 
\brkt{
\brkt{1+\eta_\ell
}^{2{|\mathbb{V}(B'')\setminus \mathbb{V}(B')|}}
\brkt{
\brkt{1+\eta_\ell
\over
1-\eta_\ell
}^{2^{|V(B')|}}+1
}
-2\brkt{1-\eta_\ell
}^{2{|\mathbb{V}(B'')\setminus \mathbb{V}(B')|}}
\brkt{1-\eta_\ell
\over
1+\eta_\ell
}^{2^{|V(B')|}}
}
\\
&\stackrel{(\ref{xx12})}{\le}&
\brkt{
\prod_{e\in\mathbb{V}(B'')\setminus\mathbb{V}(B')}{\bf d}_{\bf G}(B\ang{e})
}^2
\cdot
8.1\cdot 2^{|V(B'')|}\eta_{\ell}
\end{eqnarray*}
since \begin{eqnarray}
\eta_1\le \cdots \le
\eta_{\ell}\ll 1/k,1/\Delta.
\label{xx12}
\end{eqnarray}
It completes the claim.
\end{proof}
Fix $u\in B.$
For a set of positive integers $A$, we denote by $N_A$ the ($k$-bound)
subcomplex of $B$ induced by the 
set of vertices $v$ 
whose distances from $u$ belong to $A$ in 
the ordinary(i.e. $2$-uniform) graph $B^{\ang{2}}.$ Dropping the symbol 
$\{\}$ we simply write $N_{\{a,b\}}=N_{a,b}.$
(Note that there is no visible (hyper)edge in $B$ containing vertices from both of $N_1$ 
and $N_3$, since $B$ is a complex.)
\par
For $\ell<k$ and $i<k$, we say that $\varphi\in\Phi(N_{2i-1})$ is {\bf $\ell$-bad} iff \\
(i) $N_{2i-1}^{\ang{\ell}}\stackrel{\varphi}{\hookrightarrow}{\bf G}$ \mbox{ but, }
\\
(ii) 
 $\Ex_{\phi\in \Phi(N_{2i+1})}
[\beta(\varphi\dot{\cup}\phi,
N_{2i-1,2i+1}^{\ang{\ell}},N_{[2i-1,2i+1]}^{\ang{\ell}})|
N_{2i-1,2i+1}^{\ang{\ell}}\stackrel{\varphi\dot{\cup}\phi}{\hookrightarrow}{\bf G}
]>{\beta_\ell}^{2/3}
$.
\\
Since $|V(N_{[2i-1,2i+1]})|\le \Delta^{2k}$, 
we can apply Claim \ref{xx12b} and obtain that 
\begin{eqnarray}
&&
\Prob_{\varphi\in\Phi(N_{2i-1})}
[\varphi \mbox{ is $\ell$-bad}|N^{\ang{\ell}}_{2i-1}
\stackrel{\varphi}{\hookrightarrow}{\bf G}]
\nonumber\\
&\le &
\Prob_{\varphi
}\left[\left.
\Ex_{\phi\in \Phi(N_{2i+1})}
[\beta(\varphi\dot{\cup}\phi,
N_{2i-1,2i+1}^{\ang{\ell}},N_{[2i-1,2i+1]}^{\ang{\ell}})|
N_{2i-1,2i+1}^{\ang{\ell}}\stackrel{\varphi\dot{\cup}\phi}{\hookrightarrow}{\bf G}
]\cdot {\beta_\ell}^{-2/3}>1
\right| N^{\ang{\ell}}_{2i-1}\stackrel{\varphi}{\hookrightarrow}{\bf G}\right]
\nonumber\\
&\le &
\Ex_{\varphi
}\left[\left.
\Ex_{\phi\in \Phi(N_{2i+1})}
[\beta(\varphi\dot{\cup}\phi,
N_{2i-1,2i+1}^{\ang{\ell}},N_{[2i-1,2i+1]}^{\ang{\ell}})|
N_{2i-1,2i+1}^{\ang{\ell}}\stackrel{\varphi\dot{\cup}\phi}{\hookrightarrow}{\bf G}
]
\right| N^{\ang{\ell}}_{2i-1}\stackrel{\varphi}{\hookrightarrow}{\bf G}\right]
\cdot {\beta_\ell}^{-2/3}
\nonumber\\
&\stackrel{C.\ref{xx12b}}{\le}& \beta_\ell\cdot {\beta_\ell}^{-2/3}
=\beta_\ell^{1/3}.\label{xx14b}
\end{eqnarray}
For a $\varphi\in\Phi(N_{2i-1}),$ we define the {\bf rank} 
of $\varphi$, ${\rm rank}(\varphi)\in [0,k]$, as follows:
\\
(i) ${\rm rank}(\varphi):=0$ iff $N_{2i-1}
\stackrel{\varphi}{\hookrightarrow}{\bf G}$ does 
not hold,
\\
(ii) ${\rm rank}(\varphi):=k$ if $i=1$ and $\beta(\varphi,N_1,N_{0,1})\le 
{\beta_k}^{2/3}.$
\\
(iii) otherwise, 
 ${\rm rank}(\varphi)$ is the largest $\ell\in [k-1]$ such that 
$\varphi$ is not $\ell'$-bad for all  $\ell'\le \ell$.
\\
(Note that there is no $\varphi$ which is $1$-bad, because 
$\beta(\varphi\dot{\cup}\phi,N^{\ang{1}}_{2i-1,2i+1},N^{\ang{1}}_{[2i-1,2i+1]})=0$ 
for any $\phi\in\Phi(N_{2i+1})$.
)
\\
It follows for $1\le \ell\le k-2$ that 
\begin{eqnarray}
\Prob_{\varphi\in\Phi(N_{2i-1})}[{\rm rank}(\varphi)=\ell]
&\le& 
\Prob_{\varphi\in\Phi(N_{2i-1})}[\varphi \mbox{ is $(\ell+1)$-bad}
]
\nonumber\\
&\stackrel{(\ref{xx14b})}{\le}& 
\beta_{\ell+1}^{1/3}
\Prob_{\varphi\in\Phi(N_{2i-1})}[N_{2i-1}^{\ang{\ell+1}}
\stackrel{\varphi}{\hookrightarrow}{\bf G}].
\label{xx14}
\end{eqnarray}
A calculation similar to (\ref{xx14b}) with Claim \ref{xx12b} yields that 
\begin{eqnarray}
\Prob_{\varphi\in\Phi(N_1)}[{\rm rank}(\varphi)\in [k-1]]
\le {\beta_k^{1/3}}\Prob_{\varphi\in\Phi(N_1)}
[N_1\stackrel{\varphi}{\hookrightarrow}{\bf G}].\label{xx15c}
\end{eqnarray}
For a $\varphi\in\Phi(N_{2i-1})$ with rank $\ell\in [k-1]$ and for 
$\ell'\in [\ell]$, 
we say that $\psi\in\Phi(N_{2i+1})$ is {\bf $\ell'$-$\varphi$-bad} iff \\
(i) $N_{2i+1}^{\ang{\ell'}}\stackrel{\psi}{\hookrightarrow}{\bf G}$ 
(thus, $N_{2i-1,2i+1}^{\ang{\ell'}}
\stackrel{\varphi\dot{\cup}\psi}{\hookrightarrow}{\bf G}$)
\mbox{ but, }
\\
(ii) 
$\beta(\varphi\dot{\cup}\psi,
N_{2i-1,2i+1}^{\ang{\ell'}},N_{[2i-1,2i+1]}^{\ang{\ell'}})>
\beta_{\ell'}^{1/3}$.
\par
Furthermore for each $i$, 
the {\bf $\varphi$-rank} of $\psi\in\Phi(N_{2i+1}),$ denoted by 
${\rm rank}^{\varphi}(\psi)$, is defined 
as follows:\\
(i) the rank is $0$ iff $N_{2i+1}\stackrel{\psi}{\hookrightarrow}{\bf G}$ does not hold,\\
(ii) otherwise, it is the largest $\ell'\in [\ell]$ such that 
$\psi$ is not $\ell''$-$\varphi$-bad for all $\ell''\le \ell'.$
\\
(Note that there is no $1$-$\varphi$-bad $\psi$.)
\par
For $\ell'\in [\ell-1]$, we see that 
\begin{eqnarray}
&&
\Prob_{\psi\in\Phi(N_{2i+1})}[{\rm rank}^\varphi(\psi)=\ell']
\nonumber\\
&{\le}& 
\Prob_{\psi\in\Phi(N_{2i+1})}
[N_{2i+1}^{\ang{\ell'+1}}\stackrel{\varphi}{\hookrightarrow}{\bf G}]
\cdot
\Prob_{\psi\in\Phi(N_{2i+1})}[
\psi\mbox{ is $(\ell'+1)$-$\varphi$-bad }
|N_{2i+1}^{\ang{\ell'+1}}\stackrel{\psi}{\hookrightarrow}{\bf G}
]
\nonumber\\
&=&
\Prob_{\psi\in\Phi(N_{2i+1})}
[N_{2i+1}^{\ang{\ell'+1}}\stackrel{\varphi}{\hookrightarrow}{\bf G}]
\cdot
\Prob_{\psi\in\Phi(N_{2i+1})}[
\beta(\varphi\dot{\cup}\psi,N^{\ang{\ell'+1}}_{2i-1,2i+1},
N^{\ang{\ell'+1}}_{[2i-1,2i+1]}
)\cdot\beta_{\ell'+1}^{-1/3}>1
|N_{2i+1}^{\ang{\ell'+1}}\stackrel{\psi}{\hookrightarrow}{\bf G}
]
\nonumber\\
&\le & 
\Prob_{\psi\in\Phi(N_{2i+1})}
[N_{2i+1}^{\ang{\ell'+1}}\stackrel{\varphi}{\hookrightarrow}{\bf G}]
\cdot
\Ex_{\psi\in\Phi(N_{2i+1})}[
\beta(\varphi\dot{\cup}\psi,N^{\ang{\ell'+1}}_{2i-1,2i+1},
N^{\ang{\ell'+1}}_{[2i-1,2i+1]}
)
|N_{2i+1}^{\ang{\ell'+1}}\stackrel{\psi}{\hookrightarrow}{\bf G}
]\cdot\beta_{\ell'+1}^{-1/3}
\nonumber\\
&{\le}& 
\Prob_{\psi\in\Phi(N_{2i+1})}[N_{2i+1}^{\ang{\ell'+1}}
\stackrel{\varphi}{\hookrightarrow}{\bf G}]
\cdot
\beta_{\ell'+1}^{2/3}\cdot \beta_{\ell'+1}^{-1/3}\
\quad 
\hspace{5mm}
(\because 
\mbox{$\varphi$ is not $(\ell'+1)$-bad and $\ell'+1\le \ell<k.$})
\nonumber\\
&=& \beta_{\ell'+1}^{1/3}
\Prob_{\psi\in\Phi(N_{2i+1})}[N_{2i+1}^{\ang{\ell'+1}}
\stackrel{\varphi}{\hookrightarrow}{\bf G}].
\label{xx14a}
\end{eqnarray}
\par
For $\phi\in\Phi(B\setminus u)$ with $B\setminus u\stackrel{\phi}{\hookrightarrow}{\bf G}$,
we define the {\bf label} of $\phi$, ${\rm label}(\phi)=(a_1,\cdots,a_{k-1})\in 
[k]^{k-1}$ so that \\
$\bullet$ each $a_i$ is the minimum number among the ranks of 
$\phi|_{N_1},\phi|_{N_3},\cdots,\phi|_{N_{2i-1}}$ and among the 
$(\phi|_{N_{2j-3}})$-rank of $\phi|_{N_{2j-1}}, j=2,\cdots,i$.
\par
Since any label is a non-increasing sequence, if 
$a_1={\rm rank}(\phi|_{N_1})\le k-1$ then 
it satisfies  that 
\\
(a) 
$a_{k-1}=1$ or \\
(b) $a_{i-1}=a_i=\ell, \exists i>1,\exists \ell\in [2,k-1].$
\par
Thus we define the {\bf type} of $\phi\in\Phi(B\setminus u)$, denoted by 
${\rm type}(\phi)\in [0,k]$ as follows.\\
(i) 
The type is $0$ iff $B\setminus u\stackrel{\phi}{\hookrightarrow}{\bf G}$ does not hold.
\\
(ii) The type is $k$ if ${\rm rank}(\phi|_{N_1})=k.$\\
(iii) Otherwise,  if condition (a) holds then
 ${\rm type}(\phi):=1$, and  if (a) does not hold but 
(b)  holds then  
${\rm type}(\phi)$ is 
the largest $\ell\in [2,k-1]$ with property (b).
\begin{claim}[Case of type 1]\label{xx15f}
For any fixed $\varphi\in\Phi(N_1)$ with its rank in $[k-1]$,
 we have 
$$
\Prob_{\phi\in \Phi(N_{[2,\infty)})}[{\rm type}(\varphi\dot{\cup}\phi)=1
]\le 
\eta_2^{1\over 3.001}
\Prob_{\psi\in\Phi(B\setminus u)}[B\setminus u\stackrel{\psi}{\hookrightarrow}{\bf G}
|N_1\stackrel{\psi}{\hookrightarrow}{\bf G}
].
$$
\end{claim}
\begin{proof}
We divide it into two cases:
\\
(i) ${\rm rank}(\phi|_{N_{2i-1}})=1$ for $i\ge 1$ or \\
(ii) ${\rm rank}(\phi|_{N_{2i-3}})\ge 2$ but ${\rm rank}^{\phi|_{N_{2i-3}}}(
\phi|_{N_{2i-1}})
=1$ for some $i\ge 2.$
\\
Therefore it follows from (\ref{xx14}) and (\ref{xx14a}) 
and from $|V(N_{[2,2i-1]})|\le \Delta^{2k}$ 
that 
\begin{eqnarray*}
&&
\Prob_{\phi\in \Phi(N_{[2,\infty)})}[{\rm type}(\varphi\dot{\cup}\phi)=1
]/\Prob_{\psi\in\Phi(B\setminus u)}[B\setminus u\stackrel{\psi}{\hookrightarrow}{\bf G}
|N_1\stackrel{\psi}{\hookrightarrow}{\bf G}
]
\\
&\le &
\sum_i
{
\Prob_{\phi\in \Phi(N_{[2,2i-2]})}[N_{[2,2i-2]}^{\ang{1}}
\stackrel{\phi}{\hookrightarrow}{\bf G}
]
\Prob_{\varphi_i\in \Phi(N_{2i-1})}[\mbox{ (i) or (ii)}]
\Prob_{\psi\in \Phi(N_{[2i,\infty)})}[N_{[2i,\infty)}
\stackrel{\psi}{\hookrightarrow}{\bf G}]
\over 
\Prob_{\psi\in\Phi(B\setminus u)}[B\setminus u\stackrel{\psi}{\hookrightarrow}{\bf G}
|N_1\stackrel{\psi}{\hookrightarrow}{\bf G}
]
}
\\
&\stackrel{(\ref{xx14}),(\ref{xx14a})}\le &
\sum_i
{{\beta_2^{1/3}}
\Prob_{\phi\in \Phi(N_{[2,2i-1]})}[N_{[2,2i-1]}^{\ang{1}}
\stackrel{\psi}{\hookrightarrow}{\bf G}
]
\Prob_{\psi\in \Phi(N_{[2i,\infty)})}[N_{[2i,\infty)}
\stackrel{\psi}{\hookrightarrow}{\bf G}]
\over \Prob_{\psi\in \Phi(B\setminus u)}[B\setminus u
\stackrel{\psi}{\hookrightarrow}{\bf G}]
}
\Prob_{\psi\in \Phi(N_1)}[N_1
\stackrel{\psi}{\hookrightarrow}{\bf G}]
\\
&\stackrel{I.H.,(\ref{xx11a})}{\le} &
\sum_i
{{\beta_2^{1/3}}
\prod_{e\in\mathbb{V}_1(N_{[2,2i-1]})}
{\bf d}_{\bf G}^{(+\delta)}(B\ang{e})
\over 
(1-\eta_k^{1/4})^{|V(N_{[1,2i-1]})|}
\prod_{e\in\mathbb{V}(N_{[1,2i]})\setminus\mathbb{V}(N_{2i})}
{\bf d}_{\bf G}(B\ang{e})
}
\prod_{e\in\mathbb{V}(N_1)}{\bf d}_{\bf G}^{(+\delta)}(B\ang{e})
\\
&& (\mbox{by repeating the induction hypothesis (on $|V(B\setminus\{u\})|
>|V(N_{[2i,\infty)})|
$)
 $|V(N_{[1,2i-1]})|$ times})
\\
&\stackrel{(\ref{xx15d})}{\le} &\sum_i
{\brkt{8.1\cdot 2^{\Delta^{2k}}
\eta_2
}^{1/3}
\prod_{e\in\mathbb{V}_1(N_{[2,2i-1]})}
(1+\eta_1)
\over 0.9
\prod^{|e|\ge 2}_{e\in\mathbb{V}(N_{[1,2i]})\setminus\mathbb{V}(N_{2i})}
\rho_{|e|}
}
\prod_{e\in\mathbb{V}(N_1)}(1+\eta_{|e|})
\\
&{\le} &\eta_2^{1/3.001}
\end{eqnarray*}
since $\eta_2\ll \rho_2\le \cdots\le \rho_k,1/k,1/\Delta.$
\end{proof}
\begin{claim}[Case of full types]
\label{xx15g} Let $\ell_0\in [2,k-1].$
For any fixed $\varphi\in\Phi(N_1)$ with its rank $\ell_0$, we have 
$$
\Prob_{\phi\in \Phi(N_{[2,\infty)})}[{\rm type}(\varphi\dot{\cup}\phi)=\ell_0
]
\le 
\eta_{\ell_0+1}^{-0.01}
\Prob_{\psi\in\Phi(B\setminus u)}[B\setminus u\stackrel{\psi}{\hookrightarrow}{\bf G}
|N_1\stackrel{\psi}{\hookrightarrow}{\bf G}
].
$$
\end{claim}
\begin{proof}
We see that 
\begin{eqnarray*}
&&
\Prob_{\phi\in \Phi(N_{[2,\infty)})}[{\rm type}(\varphi\dot{\cup}\phi)=\ell_0
]/\Prob_{\psi\in\Phi(B\setminus u)}[B\setminus u\stackrel{\psi}{\hookrightarrow}{\bf G}
|N_1\stackrel{\psi}{\hookrightarrow}{\bf G}
]
\\
&\le& {
\Prob_{\psi\in\Phi(N_{3})}[{\rm rank}^{\varphi}(\psi)=\ell_0]
\Prob_{\psi\in\Phi(N_{[3,\infty)}),
\phi\in\Phi(N_2)}[N_{[2,\infty)}
\stackrel{\phi\dot{\cup}\psi}{\hookrightarrow}{\bf G}
|{\rm rank}^{\varphi}(\psi|_{N_3})=\ell_0]
\over 
\Prob_{\psi\in\Phi(B\setminus u)}[B\setminus u\stackrel{\psi}{\hookrightarrow}{\bf G}
|N_1\stackrel{\psi}{\hookrightarrow}{\bf G}
]
}
\\
&\stackrel{(\ref{xx15e})}{\le}
& (1+\beta_{\ell_0}^{1/6})\prod_{e\in \mathbb{V}_{(\ell_0)}
(N_{[1,3]})\setminus\mathbb{V}(N_{1,3})}
{\bf d}_{\bf G}(B\ang{e})\cdot 
{
\Prob_{\psi\in\Phi(N_{[3,\infty)})}[N_{[3,\infty)}
\stackrel{\psi}{\hookrightarrow}{\bf G}
|{\rm rank}^{\varphi}(\psi|_{N_3})=\ell_0]
\over 
\Prob_{\psi\in\Phi(B\setminus u)}[B\setminus u\stackrel{\psi}{\hookrightarrow}{\bf G}
|N_1\stackrel{\psi}{\hookrightarrow}{\bf G}
]
}
\\
&&\cdot 
\Prob_{\psi\in\Phi(N_{3})}[{\rm rank}^{\varphi}(\psi)=\ell_0]
\hspace{5mm}
 (\mbox{because of (\ref{xx15e}) and definiton (ii) of $\ell'$-$\varphi$-badness})
\\
&\stackrel{I.H.,(\ref{xx11a})}{\le}
& (1+\beta_{\ell_0}^{1/6})\prod_{e\in \mathbb{V}_{(\ell_0)}
(N_{[1,3]})\setminus\mathbb{V}(N_{1,3})}
{\bf d}_{\bf G}(B\ang{e})\cdot 
{
\prod_{e\in\mathbb{V}(N_1)}{\bf d}^{(+\delta)}_{\bf G}(B\ang{e})
\over (1-\eta_k^{1/4})^{|V(N_{1,2})|}
\prod_{e\in\mathbb{V}(N_{[1,3]})\setminus\mathbb{V}(N_3)}
{\bf d}_{\bf G}(B\ang{e})
}
\\
&\stackrel{(\ref{xx15d})}{\le}
& (1+\beta_{\ell_0}^{1/6})
{
\prod_{e\in\mathbb{V}(N_1)}
(1+\eta_{|e|})
\over 0.9
\prod_{e\in\mathbb{V}(N_{[1,3]})\setminus\mathbb{V}(N_3)}^{|e|>\ell_0}
\rho_{|e|}
}
\\
&\stackrel{(\ref{xx18b1})}\le & 1/\eta_{\ell_0+1}^{0.01}
\end{eqnarray*}
since 
\begin{eqnarray}
\eta_{\ell_0+1}\ll \rho_{\ell_0+1}\le \cdots \le \rho_k,1/k,1/\Delta.\label{xx18b1}
\end{eqnarray}
\end{proof}
\begin{claim}[Case of degenerate types]
\label{xx15h} Let $\ell\in [2,k-2].$
For any fixed $\varphi\in\Phi(N_1)$ with its rank in $[\ell+1,k-1]$,
 we have 
$$
\Prob_{\phi\in \Phi(N_{[2,\infty)})}[{\rm type}(\varphi\dot{\cup}\phi)=\ell
]
\le 
\eta_{\ell+1}^{1\over 3.001}
\Prob_{\psi\in\Phi(B\setminus u)}[B\setminus u\stackrel{\psi}{\hookrightarrow}{\bf G}
|N_1\stackrel{\psi}{\hookrightarrow}{\bf G}
].
$$
\end{claim}
\begin{proof}
Suppose that ${\rm label}(\varphi\dot{\cup}\phi)
=(a_1,\cdots,a_{k-1})$.
Let $i_0\in [2,k-2]$ be the smallest integer with $a_{i_0}=\ell.$
It follows from $a_{i_0+1}=a_{i_0}=\ell$ 
and from the minimality of $a_{i_0}$ 
that 
\begin{eqnarray}
\min_{j\in [2,i_0+1]}\left\{
\min\{
{\rm rank}((\varphi\dot{\cup}\phi)|_{N_{2j-3}}),
{\rm rank}^{(\varphi\dot{\cup}\phi)|_{N_{2j-3}}}
((\varphi\dot{\cup}\phi)|_{N_{2j-1}})
\}
\right\}\ge \ell,\label{xx15a}
\mbox{ and }
\\
{\rm rank}(\phi|_{N_{2{i_0}-1}})=\ell \mbox{ or }
{\rm rank}(\phi|_{N_{2i_0-3}})
>
{\rm rank}^{\phi|_{N_{2i_0-3}}}(\phi|_{N_{2i_0-1}})=\ell.\label{xx15b}
\end{eqnarray}
When $s=1,\ell+1$, let $B_s$ be the complex obtained from $B\setminus u$ by 
invisualizing all 
the edges of size at least $s$ containing a vertex of $N_{2j}, j\in [i_0].$
It follows that
\begin{eqnarray*}
&&
\Prob_{\phi\in \Phi(N_{[2,\infty)})}[{\rm type}(\varphi\dot{\cup}\phi)=\ell
]/\Prob_{\psi\in\Phi(B\setminus u)}[B\setminus u\stackrel{\psi}{\hookrightarrow}{\bf G}
|N_1\stackrel{\psi}{\hookrightarrow}{\bf G}
]
\\
&\le &\sum_{i_0}
\Prob_{\phi\in\Phi(N_{[2,\infty)})}[
B_{\ell+1}\stackrel{\varphi\dot{\cup}\phi}{\hookrightarrow}{\bf G}
\mbox{ and } (\ref{xx15a}),(\ref{xx15b})
]/
\Prob_{\psi\in\Phi(B\setminus u)}[B\setminus u\stackrel{\psi}{\hookrightarrow}{\bf G}
|N_1\stackrel{\psi}{\hookrightarrow}{\bf G}
]
\\
&\le &\sum_{i_0}
\prod_{j\in [i_0]}(1+\beta_\ell^{1/6})
\prod_{e\in\mathbb{V}_{(\ell)}(N_{[2j-1,2j+1]})\setminus\mathbb{V}(N_{2j-1,2j+1})}
{\bf d}_{\bf G}(B\ang{e}) 
\quad (\because (\ref{xx15a}),(\ref{xx15e}) )
\\
&&
\cdot 
\Prob_{\phi\in\Phi(N_{[2,\infty)})}[B_{1}\stackrel{
\varphi\dot{\cup}\phi}{\hookrightarrow}{\bf G}
\mbox{ and }  (\ref{xx15a}),(\ref{xx15b})
]
/\Prob_{\psi\in\Phi(B\setminus u)}[B\setminus u\stackrel{\psi}{\hookrightarrow}{\bf G}
|N_1\stackrel{\psi}{\hookrightarrow}{\bf G}
]
\\
&
\le &\sum_{i_0}
\prod_{j\in [i_0]}(1+\beta_\ell^{1/6})
\prod_{e}
{\bf d}_{\bf G}(B\ang{e})
\cdot 
\Prob_{\phi\in\Phi(N_{[2i_0+1,\infty)})}[N_{[2i_0+1,\infty)}
\stackrel{\phi}{\hookrightarrow}{\bf G}
]
\\
&&\cdot
\Prob_{\phi\in\Phi(N_{[2,2i_0-1]})}[B_1\setminus N_{[2i_0+1,\infty)}
\stackrel{\varphi\dot{\cup}\phi}{\hookrightarrow}{\bf G} \mbox{ and }
(\ref{xx15b})]
/\Prob_{\psi\in\Phi(B\setminus u)}[B\setminus u\stackrel{\psi}{\hookrightarrow}{\bf G}
|N_1\stackrel{\psi}{\hookrightarrow}{\bf G}
]
\\
&{\le} &
\sum_{i_0}
\prod_{j\in [i_0]}(1+\beta_\ell^{1/6})
\prod_{e}
{\bf d}_{\bf G}(B\ang{e})
\cdot {
\Prob_{\phi\in\Phi(N_1)}[N_1
\stackrel{\phi}{\hookrightarrow}{\bf G}
]\over (1-\eta_k^{1/4})^{|V(N_{[2i_0]})|}
\prod_{e\in\mathbb{V}(N_{[2i_0+1]})\setminus\mathbb{V}(N_{2i_0+1})}
{\bf d}_{\bf G}(B\ang{e})} \quad (\because \mbox{ I.H.})
\\
&&\cdot 2
\beta_{\ell+1}^{1/3}
\Prob_{\varphi\in\Phi(N_{2i_0-1})}[N_{2i_0-1}^{\ang{\ell+1}}
\stackrel{\varphi}{\hookrightarrow}
{\bf G}]
\prod_{j\in [2,i_0-1]}
\Prob_{\varphi\in\Phi(N_{2j-1})}[N_{2j-1}\stackrel{\varphi}{\hookrightarrow}
{\bf G}] \quad (\because (\ref{xx14}),(\ref{xx14a}))
\\
&\stackrel{(\ref{xx11a})}{\le} &\sum_{i_0}
\prod_{j\in [i_0]}(1+\beta_\ell^{1/6})
\prod_{e\in\mathbb{V}_{(\ell)}(N_{[2j-1,2j+1]})\setminus\mathbb{V}(N_{2j-1,2j+1})}
{\bf d}_{\bf G}(B\ang{e})
\cdot {
\prod_{e\in\mathbb{V}(N_1)}{\bf d}_{\bf G}^{(+\delta)}(B\ang{e})
\over 0.9
\prod_{e\in\mathbb{V}(N_{[2i_0+1]})\setminus\mathbb{V}(N_{2i_0+1})}
{\bf d}_{\bf G}(B\ang{e})} 
\\
&&\cdot 2
\beta_{\ell+1}^{1/3}
\prod_{e\in \mathbb{V}_{(\ell+1)}(N_{2{i_0}-1})}{\bf d}_{\bf G}^{(+\delta)}(B\ang{e})
\prod_{j\in [2,i_0-1]}
\prod_{e\in \mathbb{V}(N_{2j-1})}{\bf d}_{\bf G}^{(+\delta)}(B\ang{e})
\\
&{\le} & \sum_{i_0}
\prod_{j\in [i_0]}(1+\beta_\ell^{1/6})
\cdot {
\prod_{e\in\mathbb{V}(N_1)}(1+\eta_{|e|})
\over 0.9
\prod_{e\in\mathbb{V}(N_{[2i_0+1]})\setminus\mathbb{V}(N_{2i_0+1})}^{|e|>\ell}
\rho_{|e|}
}
\\
&&\cdot
\eta_{\ell+1}^{1/3.0001}
\prod_{e\in \mathbb{V}_{(\ell+1)}(N_{2{i_0}-1})}
(1+\eta_{|e|})
\prod_{j\in [2,i_0-1]}
\prod_{e\in \mathbb{V}(N_{2j-1})}
(1+\eta_{|e|})
\\
&\stackrel{(\ref{071122})}{\le} &\eta_{\ell+1}^{1/3.001}
\end{eqnarray*}
where we used, in the last two inequalities, the assumption that 
\begin{eqnarray}
\eta_{\ell+1}\ll \rho_{\ell+1}
\le \rho_{\ell+2}\le \cdots \le\rho_{k}
,1/k,1/\Delta.
\label{071122}
\end{eqnarray}
\end{proof}
Finally we obtain the inequalities that 
\begin{eqnarray*}
&&
\Prob_{\phi\in \Phi(B)}[B\stackrel{\phi}{\hookrightarrow}{\bf G}]
\\
&\ge&\Prob_{\varphi\in \Phi(N_1)}[{\rm rank}(\varphi)=k]
\Prob_{\phi\in\Phi(B)}[B\stackrel{\phi}{\hookrightarrow}{\bf G}|
{\rm rank}(\phi|_{N_1})=k
]
\\
&\stackrel{
(\ref{xx15e})
}{\ge}&
(1-\beta_k^{1/3})\prod_{e\in\mathbb{V}(B):u\in e}{\bf d}_{\bf G}(B\ang{e})
\cdot
\Prob_{\phi\in\Phi(B\setminus u)}[B\setminus u
\stackrel{\phi}{\hookrightarrow}{\bf G} \mbox{ and }
{\rm rank}(\phi|_{N_1})=k
]
\\
&{\ge}&
(1-\beta_k^{1/3})\prod_{e\in\mathbb{V}(B):u\in e}{\bf d}_{\bf G}(B\ang{e})
\cdot\brkt{
\Prob_{\phi\in\Phi(B\setminus u)}
[B\setminus u\stackrel{\phi}{\hookrightarrow}{\bf G}]-
\sum_{\ell_0\in [k-1]}
\Prob_{\phi\in\Phi(B\setminus u)}[{\rm type}(\phi)=\ell_0
]}
\\
&\stackrel{(\ref{xx15i})}{\ge}&
(1-\eta_k^{1/4.1})\prod_{e\in\mathbb{V}(B):u\in e}{\bf d}_{\bf G}(B\ang{e})
\cdot\brkt{1-\eta_k^{1/2.2}
}
\Prob_{\phi\in\Phi(B\setminus u)}[B\setminus u\stackrel{\phi}{\hookrightarrow}{\bf G}]
\\
&{\ge}&
\brkt{1-\eta_k^{1/4.2}}
\prod_{e\in\mathbb{V}(B):u\in e}{\bf d}_{\bf G}(B\ang{e})
\cdot
\Prob_{\phi\in\Phi(B\setminus u)}[B\setminus u\stackrel{\phi}{\hookrightarrow}{\bf G}]
\end{eqnarray*}
where we used the fact that
\begin{eqnarray}
&&
\sum_{\ell_0\in [k-1]}\Prob_{\varphi\in\Phi(N_1)}[{\rm rank}(\varphi)=\ell_0]
\nonumber\\
&&\cdot
\Prob_{\varphi\in\Phi(N_1)}[
\brkt{
\Prob_{\phi\in\Phi(N_{[2,\infty)})}
[{\rm type}(\varphi\dot{\cup}\phi)=\ell_0]
+
\sum_{\ell\in [\ell_0-1]}\Prob_{\phi
}
[{\rm type}(\varphi\dot{\cup}\phi)=\ell]
}
|{\rm rank}(\varphi)=\ell_0
]
\nonumber\\
&
\stackrel{(\ref{xx14}),(\ref{xx15c})}{\le}
& 
\sum_{\ell_0\in [k-1]}
\beta_{\ell_0+1}^{1/3}
\Prob_{\varphi\in\Phi(N_1)}[N_1^{\ang{\ell_0+1}}\stackrel{\varphi}{\hookrightarrow}
{\bf G}]
\nonumber\\
&&\cdot
\brkt{
\eta_{\ell_0+1}^{-0.01}
+
\sum_{\ell\in [\ell_0-1]}
\eta_{\ell+1}^{1/3.001}
}\Prob_{\psi\in \Phi(B\setminus u)}
[B\setminus u\stackrel{\psi}{\hookrightarrow}{\bf G}
|N_1\stackrel{\psi}{\hookrightarrow}{\bf G}
]
 \quad (\because \mbox{Claims } \ref{xx15f},\ref{xx15g},\ref{xx15h})
\nonumber\\
&\stackrel{(\ref{xx11a})}{\le}& 
\Prob_{\psi\in \Phi(B\setminus u)}
[B\setminus u\stackrel{\psi}{\hookrightarrow}{\bf G}
]
\sum_{\ell_0\in [k-1]}
\eta_{\ell_0+1}^{1\over 3.001}
\brkt{
\eta_{\ell_0+1}^{-0.01}
+
\sum_{\ell\in [\ell_0-1]}
\eta_{\ell+1}^{1\over 3.001}
}
{\prod_{e\in\mathbb{V}_{(\ell_0+1)}(N_1)}{\bf d}_{\bf G}^{(+\delta)}(B\ang{e})
\over
\prod_{e\in\mathbb{V}(N_1)}{\bf d}_{\bf G}^{(-\delta)}(B\ang{e})
}
\nonumber\\
&\le & 
\Prob_{\psi
}
[B\setminus u\stackrel{\psi}{\hookrightarrow}{\bf G}
]
\sum_{\ell_0\in [k-1]}
\eta_{\ell_0+1}^{1\over 3.001}
\brkt{
\eta_{\ell_0+1}^{-0.01}
+
\sum_{\ell\in [\ell_0-1]}
\eta_{\ell+1}^{1\over 3.001}
}
{\prod_{e\in\mathbb{V}_{(\ell_0+1)}(N_1)}(1+\eta_{|e|})
\over
\prod_{e\in\mathbb{V}(N_1)}(1-\eta_{|e|})
\prod_{e\in\mathbb{V}(N_1)}^{|e|\ge \ell_0+2}\rho_{|e|}
}
\nonumber\\
&\le & 
\Prob_{\psi\in \Phi(B\setminus u)
}
[B\setminus u\stackrel{\psi}{\hookrightarrow}{\bf G}
]
\eta_k^{1/3.1},\label{xx15i}
\end{eqnarray}
since 
\begin{eqnarray*}
\eta_{\ell_0+1}\ll \rho_{\ell_0+1}\le 
\rho_{\ell_0+2}\le \cdots \le \rho_k,
1/k,1/\Delta.
\end{eqnarray*}
It completes the proof of the main lemma.
\qed
\section{Proof of the Main Theorem}
Let $B$ be a $k$-uniform hypergraph on $n$ vertices with maximum degree 
$\Delta$, where each vertex is contained in at most ${\Delta\choose k-1}$ 
of size-$k$ visible \lq white\rq\ edges, and all non-white edges are invisible.
It is clear that $B$ can be seen as an $r$-partite hypergraph on $V(B)=V_1(B)\dot{\cup}
\cdots\dot{\cup}V_r(B)$ where $r=\Delta+1$.
(In other words, $B$ is a $\Delta$-blowup of the $k$-uniform complete hypergraph 
on $r$ vertices.)
Let ${\bf G}$ be a $k$-uniform hypergraph on $mN$ vertices, where 
each size-$k$ edge has one of $b_k$ visible colors.
Our purpose is to find a monochromatic copy of $B$ in ${\bf G}$.
We set the following parameters
$$r=\Delta+1,k,b_k \ll m\ll 1/\alpha\ll (1/\eta_i(\cdot))_{i\in [k]}
\ll 1/\varepsilon(\cdot,\cdot)
\ll \widetilde{b}_{k-1}\le \cdots\le \widetilde{b}_1
$$
with an auxiliary function 
$$\rho(\cdot)=\rho(b_i^*)=\alpha/b_i^*
$$
which will be used at (\ref{xx17a0}),(\ref{xx17a}),(\ref{xx17b}),(\ref{xx17c}),
(\ref{xx17d}).
\par 
We set $V({\bf G})={\bf \Omega}_1\dot{\cup}
\cdots\dot{\cup}{\bf \Omega}_m$ with $|{\bf \Omega}_i|=N$, and 
delete all \lq non-partitionwise\rq\ edges.
That is, any edge contains at most one vertex in a partite set ${\bf \Omega}_i$.
And color in black all the edges of size at most $k-1$.
For this resulting $m$-partite $k$-bound $(1,\cdots,1,b_k)$-colored 
graph, we apply 
 the regularity lemma (Theorem \ref{r060721}) 
with $r=m,k=k,h=2\Delta^{2k},\vec{b}=(1,\cdots,1,b_k),$ and 
\begin{eqnarray}
\varepsilon(\cdot,\cdot)=
\varepsilon(i,b^*_i)
=\alpha
\rho(b^*_i)\eta_i(\rho_i(b^*_i))
,\label{xx17a0}
\end{eqnarray}
 and obtain an 
$(\varepsilon(\cdot,\cdot),2\Delta^k)$-regular subdivision 
${\bf G}^*$ which is $(\widetilde{b}_1,\cdots,
\widetilde{b}_{k-1},\widetilde{b}_k=b_k\ge 2)$-colored where 
\begin{eqnarray}
m,k,\Delta,b_k,1/\varepsilon(\cdot,\cdot)\ll 
\widetilde{b}_1,\cdots,\widetilde{b}_{k-1}.\label{xx17a}
\end{eqnarray}
Let 
\begin{eqnarray}
b_i^*:=\max_{J:i\le |J|\le k}|{\rm C}_J({\bf G}^*)|\le \widetilde{b}_i
\mbox{ and }
\rho_i(b_i^*)= {\alpha\over b_i^*
}\ge {\alpha\over 
\widetilde{b}_i}
.\label{xx17f}
\end{eqnarray}
A size-$i$ edge ${\bf e}$ 
is called {\bf exceptional} 
iff ${\bf d}_{{\bf G}^*}({\bf G}^*\ang{\bf e})<\rho_i(b_i^*)$ or 
$\delta({\bf G}^*\ang{\bf e})>\varepsilon(i,b_i^*)/\alpha=
\rho_i(b_i^*)\eta_i(\rho_i(b_i^*))
$ where 
$\delta(\cdot)$ is a function associated with ${\bf G}^*$.
For any index $I$, it easily follows that 
$$\Prob_{{\bf e}\in {\bf \Omega}_I}[{\bf G}^*\ang{\bf e}
\mbox{ is exceptional}]\le \rho_i(b_i^*)b_i^*+
{
\varepsilon(i,b_i^*)\over 
\varepsilon(i,b_i^*)/\alpha
}=\alpha+\alpha=2\alpha.
$$
\par
Take $m$ vertices ${\bf v}_i\in {\bf \Omega}_i, i\in [m],$ randomly.
Then in the average, the number of exceptional edges of the hypergraph 
induced by the $m$ vertices 
is at most $\sum_{i\in [k]} {m\choose i}2\alpha<1$ since 
\begin{eqnarray}
m\ll 1/\alpha.\label{xx17b}
\end{eqnarray}
Thus there exist  $m$ vertices ${\bf v}_i\in {\bf \Omega}_i, i\in [m],$ 
such that 
all the edges in the graph induced by them 
are not exceptional. 
By Ramsey Theorem, Theorem \ref{xx17}, 
with 
\begin{eqnarray}
r=\Delta+1\ll m,\label{xx17c}
\end{eqnarray} 
there exist $r$ vertices among the $m$ vertices 
such that in the induced hypergraph, all of the size-$k$ edges have the same 
color, say red.
Consider $S\in {\mathcal S}_{r,k,1},$ the 
$k$-bound $r$-partite complex on {\em those} $r=\Delta+1$ vertices such that 
the color of each edge of $S$ is given by the corresponding color in 
${\bf G}^*$. (Note that all size-$k$ edges of $S$ are red.)
Denote again by $B$
 the complex obtained from the given $B$ 
\\ 
(i) by recoloring each size-$k$ white edge of $B$ in red, and\\
(ii) by coloring each edge of $B$ of size at most $k-1$ in the color of 
corresponding edge in $S$ 
so that $B$ is a $\Delta$-blowup of $S$.\par
Finally we can apply Corollary \ref{071020} (with 
$r=r,k=k,h=1,\Delta=\Delta,\vec{b}=\vec{b}^*,\rho_i=\rho_i=\alpha/b_i^*,S=S,B=B,
{\bf G}={\bf G}^*$
) where 
\begin{eqnarray}
k,\Delta,b_i^*,\cdots,b_{k-1}^*,b_k^*=b_k,1/\alpha
\ll 1/\eta_i(\rho_i(b_i^*)).
\label{xx17d}
\end{eqnarray}
We get the desired injection 
$\varphi\in\Phi(B)$ which embeds $B$ in ${\bf G},$ 
yielding a red copy of the original $k$-uniform hypergraph $B$.
\par
The above argument can be applied for any $B$ as far as 
\begin{eqnarray*}
|V_i(B)|<\eta_1(\rho_1(b_1^*))|V_i({\bf G}^*)| \quad \mbox{ for each }i
\in [r]
\end{eqnarray*}
in which by (\ref{xx17f}) the right hand side is at least 
\begin{eqnarray}
\eta_1(\rho_1(\widetilde{b}_1))N
={
\eta_1(\rho_1(\widetilde{b}_1))\over m}
{|V({\bf G})|}
=\Theta(|V({\bf G})|).
\end{eqnarray}
It completes the proof of Theorem \ref{xx16x}.
\qed
\\


\begin{thebibliography}{99}
\bibitem{BE}S.A.~Burr and P.~Erd\"{o}s, On the magnitude of generalized Ramsey numbers for 
graphs, {\em Infinite and finite sets (Colloq., Keszthely, 1973; dedicated to P.Erd\"{o}s on 
	his 60th birthday), Vol.1}, North-Holland, Amsterdam, 1975, pp.215-240. 
	Colloq.Math.Soc.J\'anos Bolyai, Vol.10.
\bibitem{Chung}F.R.K.~Chung, Quasi-random classes of hypergraphs, {\em Random Structures
Algorithms} {\bf 1}, No.4 (1990), 363-382.
\bibitem{Chung91}F.R.K.~Chung, Regularity lemmas for hypergraphs and
quasi-randomness, {\em Random Structures and Algorithms} {\bf 2} (1991), 241-252.
\bibitem{CG}F.R.K.~Chung and R.L.~Graham, Quasi-random hypergraphs,
{\em Random Structures and Algorithms} {\bf 1} No.1 (1990), 105-124.
\bibitem{CG2}F.R.K.~Chung and R.L.~Graham, Quasi-random set systems,
{\em J.~Amer.~Math.~Soc.} {\bf 4} No.1 (1991), 151-196.
\bibitem{CG3}F.R.K.~Chung and R.L.~Graham, On hypergraphs having
evenly distributed subhypergraphs,
{\em Disc. Math.} {\bf 111} (1993), 125-129.
\bibitem{CT}F.R.K.~Chung and P.~Tetali, Communication complexity and quasi randomness,
{\em SIAM J.~Discrete~Math.} {\bf 6} No.1 (1993), 110-123.
\bibitem{CRST}V. Chv\'atal, V. R\"{o}dl, E. Szemer\'edi, and
W.T.Trotter, Jr., The Ramsey number of a graph with bounded maximum
degree, {\it J.~Combin.~Theory}, {\bf B 34} (1983), 239-243.
\bibitem{ChSz}V.~Chv\'atal and E.~Szemer\'edi, 
On the Erd\H{o}s-Stone theorem, {\it J. London Math. Soc. (2)}, 207-214 (1981).
\bibitem{CFKO}O.~Cooley, N.~Fountoulakis, D.~K\"{o}hn, and D.~Osthus,
3-uniform hypergraphs of bounded degree have linear Ramsey numbers, 19 pages, 
preprint, {\tt 
arXiv:math/0608442v1} [math.CO].
\bibitem{CFKO2}O.~Cooley, N.~Fountoulakis, D.~K\"{o}hn, and D.~Osthus,
Embeddings and Ramsey numbers of sparse k-uniform hypergraphs,  
preprint, 
{\tt 
arXiv:math/0612351v1} [math.CO].
\bibitem{FR02}P.~Frankl ad V.~R\"{o}dl, Extremal problems on set systems, 
{\it Random Structures and Algorithms}, {\bf 20}(2), 131-164 (2002).
\bibitem{G}W.T.~Gowers, Hypergraph regularity and the multidimensional Szemer\'edi 
	theorem, 42 pages, preprint (2005.4, 2nd ver.)
\bibitem{GRR}R.L.~Graham, V.~R\"{o}dl and A.~Ruci\'nski, On graphs with linear 
Ramsey numbers, {\em J. Graph Theory} {\bf 35} (2000), 176-192.
\bibitem{HT}J.~Haviland and A.G.~Thomason, Pseudo-random hypergraphs,
in \lq\lq Graph Theory and Combinatorics(Cambridge, 1988)\rq\rq {\em Discrete Math.}
{\bf 75}, No.1-3 (1989), 255-278.
\bibitem{HT2}J.~Haviland and A.G.~Thomason, On testing the \lq\lq pseudo-randomness\rq\rq
of a hypergraph, {\em Discrete Math.} {\bf 103}, No.3 (1992), 321-327.
\bibitem{I}Y.~Ishigami, Proof of a conjecture of Bollob\'as and 
Kohayakawa on the Erd\H{o}s-Stone theorem, 
{\it Journal of Combinatorial Theory} {\bf B 85}, 222-254 (2002).
\bibitem{I06}Y.~Ishigami, A simple regularization of hypergraphs, 13 pages, 
preprint, 
{\tt arXiv:math/0612838v1} [math.CO].
\bibitem{I06lr}Y. Ishigami, Linear Ramsey numbers for bounded-degree hypergraphs, 
preprint, 
{\tt arXiv:math/0612601v1} [math.CO].
\bibitem{I06m}Y.~Ishigami, Removal lemma for infinitely-many forbidden hypergraphs 
and property testing, preprint, {\tt arXiv:math/0612669v1} [math.CO].
\bibitem{KNR}Y.~Kohayakawa, B.~Nagle and V.~R\"{o}dl, 
Hereditary properties of triple systems, {\em Combinatorics, Probability and Computing},
 (2003) {\bf 12}, 155-189.
\bibitem{KSSS02} J.~Koml\'os, A.~Shokoufandeh, M.~Simonovits, and 
E.~Szemer\'edi, The regularity lemma and its applications in graph theory, 
{\it Theoretical Aspects of Computer Science.} (Edited by G.B.Khosrovshahi et al.) 
 Lecture Notes in Computer Science Vol. 2292 (2002), 84-112.
\bibitem{KR}A.V.~Kostochka and V.~R\"{o}dl, On Ramsey numbers of 
uniform hypergraphs with given maximum degree, 
{\em Journal of Combinatorial Theory}, {\bf A 113}(2006) 1555-1564.
\bibitem{NORS}B.~Nagle, S.~Olsen, V.~R\"{o}dl and M.~Schacht,
On the Ramsey number of sparse 3-graphs, 20 pages, preprint (2006).
\bibitem{NRS}B.~Nagle, V.~R\"{o}dl and M.~Schacht, 
The counting lemma for regular $k$-uniform hypergraphs, 
{\em Random Structures and Algorithms}, {\bf 28} (2006), no.2, 113-179.
\bibitem{Ramsey}F.P.~Ramsey, On a problem of formal logic,
{\it Proc. Lond. Math. Soc.} (2) {\bf 30} (1930), 264-286.
\bibitem{RS}V.R\"{o}dl and M.Schacht, Regular partitions of hypergraphs, 
{\em Combinatorics, Probability \& Computing}, to appear (preprint 50 pages, 2006.5).
\bibitem{RSk04}V.~R\"{o}dl and J.~Skokan, Regularity lemma for 
$k$-uniform hypergraphs, {\em Random Structures and Algorithms} {\bf 25} (2004) (1),
1-42.
\bibitem{Sz}
	E.~Szemer\'edi,
	{Regular partitions of graphs}
	in {\it Probl\`emes combinatoires et th\'eorie des graphes},
	Orsay 1976,
	J.-C. Bermond, J.-C. Fournier, M. Las Vergnas, D. Sotteau, eds.,
	Colloq. Internat. CNRS 260,
	Paris, 1978, 399--401.
\bibitem{Tao06}T.~Tao, A variant of the hypergraph removal lemma, 
{\em J.~Combin.~Theory} {\bf A 113} (2006), no.7, 1257-1280.
\end{thebibliography}
\end{document}